\documentclass{article}

\usepackage{amsmath}
\usepackage{amsthm}
\usepackage{amssymb}

\title{Disturbing the Dyson Conjecture\\ (in a GOOD Way)}
\author{Andrew V. Sills and Doron Zeilberger\footnote{Partially
supported by NSF grant DMS 0401124}\\
\small{Department of Mathematics, Rutgers University}\\
\small{Hill Center, Busch Campus, Piscataway, NJ 08854-8019}\\
\\ \small{2000 AMS Subject Classification 05-04}
\\ \small{KEY WORDS: Dyson's conjecture, constant term} }

\date{\small{submitted August 2, 2005; minor revisions November 22, 2005}}
\newcommand{\CT}{\mathrm{CT}}
\newcommand{\xv}{\mathbf{x}}
\newcommand{\xvk}{\mathbf{\hat{x}}^{(k)}}
\newcommand{\av}{\mathbf{a}}
\newcommand{\avk}{\mathbf{\hat{a}}^{(k)}}
\newcommand{\bv}{\mathbf{b}}

\newcommand{\zerov}{\mathbf{0}}
\newcommand{\Comp}{\mathrm{Comp}}

\theoremstyle{plain}
\newtheorem*{DC}{Dyson's conjecture}

\numberwithin{equation}{section}

\begin{document}
\maketitle
\begin{center}\emph{Dedicated to Freeman John Dyson and Irving John Good}
\end{center}

\begin{abstract}
We present a case study in {\it experimental} yet {\it rigorous}
mathematics by describing an algorithm, fully implemented in
both Mathematica and Maple,
that {\it automatically conjectures}, and then
{\it automatically proves}, closed-form expressions extending Dyson's
celebrated constant term conjecture.
\end{abstract}

\section{Introduction}
Let
\begin{equation*}\label{Fdef}
F_n ( \xv; \av; \bv) := \prod_{h=1}^n x_h^{-b_h}\prod_{1\leqq i\neq j\leqq n}
  \left( 1-\frac{x_i}{x_j} \right)^{a_j}
\end{equation*}
where $\xv=\langle x_1, x_2, \dots, x_n \rangle$,
$\av=\langle a_1, a_2, \dots, a_n \rangle$,
$\bv=\langle b_1, b_2, \dots, b_n \rangle$
and
\begin{equation*}
c_n (\av, \bv) = \CT \big( F_n(\xv; \av; \bv)\big),
\end{equation*}
where $\CT (X)$ denotes the constant term, i.e. the coefficient of
$x_1^0 x_2^0 \cdots x_n^0$, in the expression $X$.
The following conjecture is due to
Freeman Dyson~\cite[p. 152, Conjecture C]{fjd}:
\begin{DC} For positive integers $n$ and nonnegative integers $a_i$,
$1\leqq i \leqq n$,
\[ \CT\left( F_n(\xv; \av; \zerov) \right)
 = \frac{ (a_1 + a_2 + \cdots + a_n)!}{a_1! a_2! \cdots a_n!}, \]
\end{DC}
\noindent where $\zerov = \langle 0,0,\dots,0 \rangle$.
Dyson noted that the $n=1, 2$ cases are trivial and that the $n=3$ case is
equivalent
to a hypergeometric summation formula due to A.C. Dixon~\cite{acd}.
Dyson proved the $n=4$ case~\cite[pp. 155--156, Appendix B]{fjd},
and indicated that a similar argument could be used to prove $n=5$, but
that his argument would not work for $n>5$, and accordingly
left $n>5$ as a conjecture.
The conjecture was
quickly proved independently by J.~Gunson~\cite{jg} and
K.~Wilson~\cite{kw}.
The most compact and elegant proof, however, was supplied by
I.~J.~Good~\cite{ijg}.
A combinatorial proof
was later given by Zeilberger~\cite{dz}.

 In this paper, we concern ourselves with variations on Dyson's original conjecture
where $\bv$ can assume (fixed) values other than $\zerov$.
In particular, we have automated, in the accompanying
Mathematica package \texttt{GoodDyson.m} and
analogous Maple package \texttt{GoodDyson},
both the act of \emph{conjecturing}
an explicit formula for $c_n (\av; \bv)$ (where $\bv$ is a vector of specific
integers), and the production of a \emph{proof} of the conjectured form,
based on a generalization
of Good's ideas.  The accompanying Mathematica and
Maple packages are
available from the authors' web sites
\texttt{http://www.math.rutgers.edu/\~{}asills} and
\texttt{http://www.math.rutgers.edu/\~{}zeilberg}.

\section{Automating the Conjecturing Process}
Given $\bv=\langle b_1, \dots, b_n \rangle$, we guess
that the coefficient of
$x_1^{b_1} x_2^{b_2} \cdots x_n^{b_n}$ in $F_n(\xv; \av; \zerov)$, or
equivalently, the constant term of $F_n(\xv; \av; \bv)$, which we denote by
$c_n (\av; \bv)$, can be expressed in the form $d_n(\av ; \bv)$, where
   \[ d_n(\av;\bv) = R_{\bv}(\av)
   \frac{(a_1 + a_2 + \cdots + a_n)!}{a_1! a_2! \cdots a_n!}, \]
and $R_{\bv}(\av)$ is a rational function in the $a_i$'s.  Of course, if
$\sum_{i=1}^n b_i \neq 0$, then $R_{\bv}(\av) = 0$, and in the case of
Dyson's original conjecture, we have $R_{\zerov} (\av) = 1$ for all $n$.

 We programmed a Mathematica function/Maple procedure
\texttt{GuessRat}, which takes
as input a function $f$, a set of variables
$\{ a_1, \dots a_n \}$, and an integer $t$,
and tries to match $f$ to a rational function $R$ in which the sum of
the degrees of the numerator and denominator (let us call this the
\emph{total degree of $R$}) is $t$, using internally
generated data.  Of course, if no such $R$ of total degree $t$ is found, then
the procedure is repeated with a larger $t$.  \texttt{GuessRat}
could potentially be useful in a wide variety of settings, but for this project
we restricted $f$ to the function which extracted the constant term
from $F(\xv; \av; \bv)/\frac{(a_1+\cdots+a_n)!}{a_1! \cdots a_n!}$
for a specific $\bv$ in order to conjecture $R_{\bv}(\av)$.  This is done
via the auxiliary function \texttt{GuessDysonCoeff}.

  Empirical evidence gathered while testing a prototype version of
\texttt{GuessDysonCoeff}
suggested that for each $b_i < 0$, the $R_{\bv}(\av)$ contained a factor of
\[ \frac{ 1 }
 { (1+a_i)_{\lfloor b_i/2\rfloor}
 (1+ a_1 + a_2 + \cdots + a_{i-1} + a_{i+1} + \cdots + a_n)_{|b_i|} },\]
where $(y)_h$ denotes the rising factorial, which is defined
by \[ (y)_h := \left\{
    \begin{array}{ll}
      y(y+1)(y+2)\cdots(y+h-1), &\mbox{if $h>0$,}\\
      1                       , &\mbox{if $h=0$,}\\
      \frac {1}{(y-1)(y-2)\cdots(y-h)}, &\mbox{if $h<0$.}
    \end{array} \right. \]

  For example, closed form representation of
the constant term of \mbox{$F_4 (\xv; \av; \langle -3, b_2, -1, 4-b_2 \rangle)$}
contains the factor
\[ \frac{a_1(a_1-1) a_3}
{(1+a_2+a_3+a_4)(2+a_2+a_3+a_4)(3+a_2+a_3+a_4)(1+a_1+a_2+a_4)} .\]

 Accordingly, we modified the $f$ which was sent as input into \texttt{GuessRat}
via the \texttt{GuessDysonCoeff} function,
resulting in the output $R$ being of lower total degree, and therefore greatly reducing
the time Mathematica/Maple needs to supply a conjecture.  For a vector $\bv$ whose
components sum to zero,
let us define the \emph{complexity} of $\bv$, $\Comp(\bv)$, to be the sum of
its positive components, or equivalently,
\[ \Comp(\bv) := \frac 12\sum_{i=1}^n |b_i|. \]
When the complexity of of $\bv$ is
close to zero, the modified algorithm
worked over twenty times faster than the original.  For larger
complexity, the speedup was even more significant.  For instance,
in the case $\langle b_1, b_2, b_3 \rangle = \langle 4,-2,-2 \rangle$, the
modified algorithm was over one hundred times faster than the original.

\section{The Generalized Good Proof}
For some fixed $n$ and $\bv$, we want to show that for all $\av$,
\begin{equation}\label{ceqd}
 c_n(\av; \bv) = d_n (\av; \bv),
 \end{equation}
where an explicit closed form expression $d_n(\av; \bv)$ for $c_n(\av; \bv)$ has been
conjectured using Mathematica or Maple.
Since the $n=1$ case is trivial, and the $n=2$ case
\begin{equation} \label{case2}
c_2(\langle a_1, a_2 \rangle; \langle b_1, b_2 \rangle)
= \left\{ \begin{array}{ll}
             \displaystyle{(-1)^{b_1}
              \frac{(a_1+a_2)!}{(a_1+b_1)!(a_2-b_1)!},} &\mbox{if $b_1=-b_2$}\\
           0, &\mbox{otherwise}
                        \end{array}
      \right.
\end{equation}
follows from the binomial theorem,
 we restrict our attention to $n>2$.

 As in Good's proof~\cite{ijg}, for $a_1, a_2, \dots, a_n \geqq 1$,
$F_n ( \xv; \av; \bv)$ satisfies the recursion
\begin{equation}\label{Frec}
  F_n( \xv; \av; \bv) = \sum_{i=1}^n
    F_n(\xv; \langle a_1, \dots, a_{i-1}, a_i - 1, a_{i+1},\dots , a_n \rangle; \bv).
\end{equation}

Applying the constant term operator to both sides of \eqref{Frec}, we obtain
\begin{equation}\label{crec}
  c_n (\av; \bv) = \sum_{i=1}^n
  c_n (\langle a_1, \dots, a_{i-1}, a_i - 1, a_{i+1},\dots , a_n \rangle; \bv).
\end{equation}

Next, we note that for any fixed $k$, $1\leqq k \leqq n$,
\begin{equation}  \label{GoodKernel}
 F_n ( \xv; \langle a_1, \dots, a_{k-1}, 0, a_{k+1}, \dots, a_n \rangle; \bv)\\
= F_{n-1} ( \xvk, \avk, \zerov)
\left[ x_k^{-b_k}\underset{i\neq k}{\prod_{i=1}^n}
  \frac{ (x_i - x_k)^{a_i} }{ x_i^{a_i+b_i} }\right],
 \end{equation}
 where
 $\xvk = \langle x_1,\dots, x_{k-1}, x_{k+1}, \dots x_n \rangle$ and
 $\avk = \langle a_1,\dots, a_{k-1}, a_{k+1}, \dots a_n \rangle$.
 Notice that on the right hand side of \eqref{GoodKernel},
 we have managed to segregate the factors
 involving $x_k$ (those in brackets) from those which do not
 involve $x_k$.  We can therefore let
Mathematica or Maple find the explicit Taylor series expansion of
$\underset{i\neq k}{\prod_{i=1}^n} \frac{ (x_i - x_k)^{a_i} }{ x_i^{a_i+b_i} }$
about $x_k = 0$.  And so, by extracting the coefficient of $x_k^0$ on both
sides of \eqref{GoodKernel}, we have, for $1\leqq k \leqq n$,
\begin{gather} \label{GK2}\mbox{coeff of $x_k^0$ in\ }
\left[ F_n ( \xv; \langle a_1, \dots, a_{k-1}, 0, a_{k+1}, \dots, a_n \rangle; \bv)
\right]\\
= \mbox{coeff of $x_k^0$ in\ }
\left[ \left( x_k^{-b_k}\underset{i\neq k}{\prod_{i=1}^n}
  \frac{ (x_i - x_k)^{a_i} }{ x_i^{a_i+b_i} } \right) F_{n-1} ( \xvk, \avk, \zerov)
  \right]
. \nonumber
\end{gather}
Finally, we apply the constant term operator to both sides of \eqref{GK2} to
obtain
\begin{gather}\label{BoundCond}
 c_n (\langle a_1, \dots, a_{k-1}, 0, a_{k+1}, \dots, a_n \rangle; \bv)\\
= P_k\ c_{n-1} ( \avk, \zerov),\nonumber
\end{gather}
where $P_k$ is
the coefficient of $x_k^{b_k}$ in the Taylor expansion of
$\underset{i\neq k}{\prod_{i=1}^n} \frac{ (x_i - x_k)^{a_i} }{ x_i^{a_i+b_i} }$
about $x_k = 0$.  Notice that $P_k$ is a Laurent polynomial in
$x_1, \dots, x_{k-1}, x_{k+1}, \dots, x_n$, whose coefficients depend on
$a_1, \dots, a_{k-1}, a_{k+1}, \dots, a_n$.
(In the case where $\bv = \zerov$, we have $P_k =1$,
making Good's proof very tidy indeed!)
Finally, we note the initial condition
\begin{equation}\label{IC}
c_n ( \zerov; \bv) = \left\{ \begin{array}{ll}
                                               1, &\mbox{if $\bv=\zerov$}\\
                                               0, &\mbox{otherwise.}
                                              \end{array}
                               \right.
\end{equation}
The equation \eqref{crec}, the set of $n$ boundary conditions \eqref{BoundCond},
together with \eqref{IC} fully determine $c_n(\av; \bv)$.  Thus, to prove \eqref{ceqd},
it suffices to show that our conjectured formula
$d_n(\av;\bv)$ also satisfies \eqref{crec}, \eqref{BoundCond},
and \eqref{IC}.  The equations \eqref{BoundCond} will, in general, depend on
$c_{n-1} (\av; \bv)$ for various values of $\bv$, and so the boundary conditions
will need to be iterated $n-2$ times until \eqref{BoundCond} is expressed fully
in terms of $c_2$'s which is given by \eqref{case2}.
\section{Example}
The coefficient of $\displaystyle{\frac{x_1^2}{x_2 x_3}}$ in the expansion of
the Laurent polynomial
\begin{gather*}
 \left[ \left(1-\frac{x_2}{x_1}\right)
          \left( 1-\frac{x_3}{x_1}\right) \right]^{a_1}
\left[ \left(1-\frac{x_1}{x_2}\right)
          \left( 1-\frac{x_3}{x_2}\right) \right]^{a_2}
\left[ \left(1-\frac{x_1}{x_3}\right)
           \left( 1-\frac{x_2}{x_3}\right) \right]^{a_3}
 \end{gather*}
    is
\begin{gather*} d_3 ( \langle a_1, a_2, a_3\rangle; \langle 2,-1,-1 \rangle )=
\frac{a_2 a_3 (2+2a_1+a_2+a_3)(a_1+a_2+a_3)!}
{(1+a_1+a_2)(1+a_1+a_3)(1+a_1) a_1! a_2! a_3! }.
\end{gather*}
\begin{proof}[Good style proof]
It is easily verified that for $a_1, a_2, a_3\geqq 1$,
\begin{eqnarray}\label{Frecex}
  F_3( \xv; \av; \langle 2,-1,-1\rangle) &=&
  F_3( \xv; \langle a_1-1, a_2, a_3 \rangle; \langle 2,-1,-1\rangle)\\
&&\qquad +F_3( \xv; \langle a_1, a_2-1, a_3 \rangle; \langle 2,-1,-1\rangle)
\nonumber\\
&&\qquad +F_3( \xv; \langle a_1, a_2, a_3-1 \rangle; \langle 2,-1,-1\rangle)
\nonumber
\end{eqnarray}
Applying the constant term operator to both sides of \eqref{Frec}, we
immediately obtain
\begin{eqnarray}\label{crecex}
  c_3(\av; \langle 2,-1,-1\rangle) &=&
  c_3(\langle a_1-1, a_2, a_3 \rangle; \langle 2,-1,-1\rangle)\\
&&\qquad +c_3( \langle a_1, a_2-1, a_3 \rangle; \langle 2,-1,-1\rangle)
\nonumber\\
&&\qquad +c_3(\langle a_1, a_2, a_3-1 \rangle; \langle 2,-1,-1\rangle)\nonumber
\end{eqnarray}

The boundary conditions are found by \texttt{GoodDyson} to be
\begin{eqnarray}
c_3(\langle 0, a_2, a_3 \rangle; \langle 2,-1,-1 \rangle) &=&
\frac{a_3(a_3-1)}{2} c_2(\langle a_2, a_3 \rangle; \langle -1,1 \rangle)\nonumber\\
&&\quad+\frac{a_2(a_2-1)}{2} c_2(\langle a_2, a_3 \rangle; \langle 1,-1 \rangle)
\nonumber\\ &&\quad
+ a_2 a_3 c_2 (\langle a_2, a_3 \rangle; \langle 0,0 \rangle)\\
c_3(\langle a_1, 0, a_3 \rangle; \langle 2,-1,-1 \rangle) &=& 0\\
c_3(\langle a_1, a_2, 0 \rangle; \langle 2,-1,-1 \rangle) &=& 0
\end{eqnarray}
Finally,
\begin{equation}\label{ICex}
c_3(\langle 0,0, 0 \rangle; \langle 2,-1,-1 \rangle) =0
\end{equation}

Since  $c_3(\av; \langle 2,-1,-1\rangle)$ is uniquely determined by
\eqref{crecex}--\eqref{ICex}, and it is easily verified that,
in light of \eqref{case2}, our conjectured expression
$d_3(\av;\langle 2,-1,-1 \rangle)$ also satisfies
\eqref{crecex}--\eqref{ICex}, the proof is complete.
\end{proof}
 Full proofs, analogous to the above, can be generated automatically with our
Mathematica function/Maple procedure \texttt{WritePaper}.

\section{Exploiting Symmetry and Algebraic Relations}
Suppose $d_n (\av; \bv)$ is known for some particular $\bv$.  It is then a simple
matter to determine $d_n (\av; \bv')$ for all vectors $\bv'$ whose components are
a permutation of the components of $\bv$: if $\pi_\bv$ permutes the indices of
$\bv=\langle b_1, b_2, \dots, b_n \rangle$, so that $\pi_\bv \bv = \bv'$, then
$d_n(\av;\bv') = d_n( \pi_\av \av; \bv)$.  For example, given
 \[ d_3(\av; \langle -1, 0, 1 \rangle) =
 -\frac{(a_1+a_2+a_3)!a_1}{a_1! a_2! a_3! (1+a_2+a_3)}, \]
 we immediately know that
  \begin{gather*} d_3(\av; \langle 1, 0, -1 \rangle) =
 -\frac{(a_1+a_2+a_3)!a_3}{a_1! a_2! a_3! (1+a_2+a_1)}, \\
   d_3(\av; \langle 0,- 1, 1 \rangle) =
 -\frac{(a_1+a_2+a_3)!a_2}{a_1! a_2! a_3! (1+a_1+a_3)},
 \end{gather*}
 etc.

  Also, notice that
 \[ F_n(\xv;\av+\langle 0,0,\dots,0,1 \rangle; \bv) =
     F_n(\xv;\av;\bv) \left[ \left( 1-\frac{x_1}{x_n} \right)
     \left( 1-\frac{x_2}{x_n} \right)\cdots
     \left( 1-\frac{x_{n-1}}{x_n} \right) \right]. \]
By expanding out the Laurent polynomial in brackets, distributing it over the
$F_n(\xv;\av;\bv)$ and applying the constant term operator, one can often
express a $d_n(\xv;\av;\bv)$ with a $\bv$ of higher complexity as a linear
combination of $d_n(\xv;\av;\bv)$'s with $\bv$'s that have previously been
calculated and/or have lower complexity.  For example,
\begin{eqnarray*}&& F_3 (\xv; \av+\langle 0,0,1 \rangle; \bv) \\
& = & F_3(\xv;\av; \bv) \left( 1-\frac{x_1}{x_3} \right)
 \left( 1-\frac{x_2}{x_3} \right) \\
 &=& F_3(\xv;\av; \bv) \left( 1-\frac{x_1}{x_3} - \frac{x_2}{x_3}
 + \frac{x_1 x_2}{x_3^2} \right)\\
 &=& F_3(\xv;\av;\bv) -F_3 (\xv;\av; \bv+\langle -1,0,1 \rangle)
 -F_3 (\xv;\av; \bv+\langle 0,-1,1 \rangle)
 +F_3 (\xv;\av; \bv+\langle -1,-1,2 \rangle)
 \end{eqnarray*}
After applying the constant term operator to both sides and solving
 for the last term, we find that
 \[ d_3(\av;\bv+\langle -1,-1,2 \rangle) =
 d_3(\langle a_1, a_2, 1+a_3 \rangle;\bv)
 - d_3(\av;\bv) + d_3 (\av; \bv+\langle -1,0,1 \rangle)
 + d_3 (\av; \bv+\langle 0,-1,1 \rangle) \]

  By systematically taking advantage of the above observations, the procedure
\texttt{TurboDyson} can find $d_n(\av;\bv)$ for a given $n$ and all $\bv$'s
of complexity less than or equal to a given complexity $C$ rather quickly.  Furthermore,
\texttt{TurboDyson} stores its findings in an indexed global variable for future
use, e.g. a subsequent call to the \texttt{WritePaper} procedure.

\section{Conclusion}
Experimental mathematics, as it is commonly understood, consists of
computer-assisted {\it gathering of data} (usually numeric, but recently
also symbolic) that stimulates and inspires {\it human-made} conjectures,
that, in turn, require {\it human-made} {\it proofs}. The novelty of
the present research is that all these phases are done by machine,
once the initial effort of designing an algorithm, and implementing it,
are done by humans. Of course, at present, the general framework
for the conjecture, and the {\it idea} of proof, are still
human-made (in our case by two very illustrious humans: Dyson and Good),
but all the rest is automatic. We believe that this methodology should
be applicable to many other problems.

\end{document}